\documentclass[12pt]{amsart}
\usepackage{xypic}
\usepackage{amsmath}
\usepackage{amssymb}
\usepackage{stmaryrd}
\newtheorem{theorem}{Theorem}[section]
\newtheorem{corollary}[theorem]{Corollary}

\theoremstyle{definition}

\theoremstyle{remark}
\newtheorem{remark}[theorem]{Remark}

\title[On the verification of a Nicolas inequality]
{On the verification of a Nicolas inequality}

\begin{document}
%----------Author 1
\author[O. Galdames-Bravo]{Orlando Galdames-Bravo}

% \address{Departament de Matem\`atiques\\
% CIPFP Vicente Blasco Ib\'a\~nez\\
% Avda. del Regne de Val\`encia, 46\\
% 46005-Val\`encia. Spain.}

\email{galdames@uv.es}

\thanks{}
%----------Author 2
% \author{A Second Author}
% \address{The address of\br
% the second author\br
% sitting somewhere\br
% in the world}
% \email{dont@know.who.knows}
%\subjclass[2020]{Primary 26D07, Secondary 11M26}


%%% ----------------------------------------------------------------------
\maketitle
%%% ----------------------------------------------------------------------
\begin{abstract}
Nicolas inequality we deal can be written as
\begin{equation}\label{Nicineq}
e^\gamma \log\log N_x < \dfrac{N_x}{\varphi(N_x)}\,,
\end{equation}
where $x\ge 2$, $N_x$ denotes the product 
of the primes less or equal than $x$, $\gamma$ is the Euler constant and 
$\varphi$ is the Euler totient function. We see that verification of  
\eqref{Nicineq} depends on the sign of the big-O function in the Mertens 
estimate for the sum of reciprocals of primes. Then, we analyze the sign 
of such an error term.
\end{abstract}

%%%%%%%%%%%%%%%%%%%%%%%%%%%%%%%%%%%%%%%%%%%%%%%%%%%%
\subsection*{Author's note}
The sequence of versions of this paper that we have reviewed illustrates the evolution 
from an astonishing result to another one that could be a good tool for studying Nicolas 
inequality. We have been able to trace the errors, mistakes, typos and inaccuracies 
from a first draft to a final version. Despite this does not guarantee the 
paper is entirely free of oversights. 

%%%%%%%%%%%%%%%%%%%%%%%%%%%%%%%%%%%%%%%%%%%%%%%%%%%%
\section{Introduction and notation}
%%%%%%%%%%%%%%%%%%%%%%%%%%%%%%%%%%%%%%%%%%%%%%%%%%%%

Nicolas inequality was introduced in \cite[Theorem 2]{Nic} and establish that 
one and only one of the following holds:\\

(T) Inequality \eqref{Nicineq} holds for every integer $x\ge 2$.\\

(F) Inequality \eqref{Nicineq} does not hold for an infinite set of integers $x$'s.\\

However, as far as we know, it is unsolved which of them is true.\\

Mertens proved in \cite{Mert} that for $x\ge 2$ 
$$
\sum_{p\le x} \frac1{p} = \log\log x + b + O(1/\log x)\,,
$$
where $b$ is sometimes called Mertens constant.\\

It turns out that if we apply adequately this Mertens estimate to the Nicola's 
inequality we see that the sign of the error term $O(1/\log x)$ is decisive 
for the compliance or the failure of Nicolas inequality. In fact, 
we prove that if $O(1/\log x)$ is bounded from below by 
a positive number, then (T) holds, unlike if $O(1/\log x)$ is bounded from above by
a negative number, then (F) holds. After that, follow the proof of Mertens Theorem
and find that we cannot determine the sign of $O(1/\log x)$ with the known estimates.  
For the sake of simplicity we use the notation $N_x=\prod_{p\le x}p$, 
instead of the original $N_k$, which is the product of 
the first $k$ primes. Notice that both are equivalent, 
because $N_x=N_{[x]}$, where $[x]$ is the floor function, 
and for any positive integer $k$ there is a real number $x$ such that 
$k=[x]$.\\

Let us introduce some terminology.
We have used standard notation for number theory that can be found, 
for instance, in \cite{Niv,Ten-Men}. 
Given a positive function $g(x)$, the \emph{big-O function} of $g(x)$,
denoted by $O(g(x))$, is a function such that $|O(g(x))|\le Cg(x)$ uniformly. Observe
that necessarily $C>0$, otherwise $O(g(x))$ vanishes, but $O(g(x))$ could be negative, 
positive or both. Moreover we can write $O(g(x))=C(x)g(x)$ for some uniformly bounded
function $C(x)$. Throughout the paper we consider $x$ as a positive
real variable, $p$ denotes prime numbers and $n$ and $d$ denote natural numbers. 
\emph{Prime zeta function} is defined as the 
sum of the powers of reciprocals of prime numbers, i.e. 
$P(s)=\sum_{p}\frac1{p^s}$, where $s$ is 
a complex number, (see for instance \cite{Fro,Mer}), 
also $P_x(s)=\sum_{p\le x}\frac1{p^s}$.
The \emph{von Mangoldt function} $\Lambda(n)$ is the arithmetic function
$\Lambda(n)=\log p$ if $n=p^k$ for some positive integer $k$, $\Lambda(n)=0$ otherwise. The
so-called \emph{Chebyshev functions} are $\psi(x)=\sum_{n\le x}\Lambda(n)$ and 
$\theta(x)=\sum_{p\le x}\log p$, and the \emph{prime counting function}
is $\pi(x)=\sum_{p\le x}1$. We define $N_x=\prod_{p\le x}p$ and
$\varphi_x=\prod_{p\le x}(p-1)$. Observe that $\theta(x)=\log N_x$ 
and $\varphi_x=\varphi(N_x)$, where 
$\varphi(n)=|\{m\in\mathbb{N}:m\le n, (m,n)=1\}|$ is the
Euler totient function. The \emph{floor function} is 
$[x]=\max\{m\in\mathbb{Z}:m\le x\}$, so we can write $[x]=x-\{x\}$, 
where function $\{x\}$ is bounded between $0$ and $1$.

%%%%%%%%%%%%%%%%%%%%%%%%%%%%%%%%%%%%%%%%%%%%%%%%%%%%
\section{On the sign of the error term}
%%%%%%%%%%%%%%%%%%%%%%%%%%%%%%%%%%%%%%%%%%%%%%%%%%%%
In this unique section we prove a theorem and discuss the error term. The first one
relates the sign of the error term of Mertens theorem for 
the sum of reciprocal of primes with the verification of the Nicolas inequality.
To this aim we apply results involving the zeta prime function,
basic estimates for Chevyshev functions and standard techniques
of real analysis. In the second part we analyze the sign
of such an error term, to this end we revise the proof of 
Mertens Theorem and obtain an explicit formula for the term,
then we approximate sufficiently its value. Unfortunately,
the conditions we impose never hold by virtue of \cite[Theorem 1.1]{But}, 
since the error term change of sign infinitely often, see also \cite{Rob}.
After that, we reformulate the theorem.

\begin{theorem}
Suppose that $x\ge 2$, and let the Mertens formula
$$
\sum_{p\le x} \frac1{p} = \log\log x + b + O(1/\log x)\,.
$$
(I) If there are $x_0\ge 2$ and $C>0$ such that $O(1/\log x)>\frac{C}{\log x}$ for every $x\ge x_0$,
then (T) holds.\\
(II) If there are $x_0\ge 2$ and $C\le0$ such that $O(1/\log x)\le\frac{C}{\log x}$ for every $x\ge x_0$,
then (F) holds.
\end{theorem}
\begin{proof} We divide the proof in three steps.\\

{\it Step 1}. With the Taylor series of the logarithm we get
\begin{equation*}
\begin{split}
\log \varphi_x & = \log\prod_{p\le x}(p-1) = 
\log\Big(\prod_{p\le x}p\prod_{p\le x}\big(1-\frac1{p}\big)\Big) \\
&= \log\prod_{p\le x}p + \log\prod_{p\le x}\big(1-\frac1{p}\big) 
= \theta(x) - \sum_{p\le x}\sum_{n=1}^\infty \frac1{np^n} \\
&= \theta(x) - \sum_{n=1}^\infty \frac1{n}\sum_{p\le x}\frac1{p^n}
= \theta(x) - \sum_{n=1}^\infty \frac{P_x(n)}{n}\,.
\end{split}
\end{equation*}

Taking into account Mertens Theorem (see e.g. \cite[Theorem 8.8(d)]{Niv})
and the above equation, we can rewrite the Nicolas inequality, namely 
$e^\gamma \log\log N_x < \dfrac{N_x}{\varphi_x}$, as
\begin{equation}\label{gamma}
\begin{split}
\gamma &< \log\Big(\frac{N_x}{\varphi_x\log\log N_x}\Big) 
=\log N_x - \log\varphi_x - \log\log\log N_x \\
&= \sum_{n=1}^\infty \frac{P_x(n)}{n} - \log\log\theta(x)
= \sum_{n=2}^\infty \frac{P_x(n)}{n} - \log\log\theta(x) + \sum_{p\le x}\frac1{p}\\
&= \sum_{n=2}^\infty \frac{P_x(n)}{n} - \log\log\theta(x) + \log\log x + b + O(1/\log x)\\
&= \sum_{n=2}^\infty \frac{P_x(n)}{n} + b 
+ \log\Big(\frac{\log x}{\log\theta(x)}\Big) + O(1/\log x)\,.
\end{split}
\end{equation}

Observe that $P_x(n)$ is increasing in $x$ and converges to $P(n)$, thus for 
each $\varepsilon>0$ there exists $x_1$ such that 
\begin{equation}\label{part}
\sum_{n=2}^\infty \dfrac{P(n)}{n}-\sum_{n=2}^\infty \dfrac{P_x(n)}{n}<\varepsilon\,,\qquad
\mbox{ for every } x\ge x_1\,.
\end{equation}
In addition one can deduce that 
\begin{equation}\label{zeta}
\sum_{n=2}^\infty \dfrac{P(n)}{n}=\gamma-b\,,
\end{equation} 
where $\gamma$ is the Euler constant and $b$ is known as the Mertens constant 
(\cite[(2.7)]{Ros-Sch}).\\

{\it Step 2}: (I). 
Let $\alpha>0$, as $C>0$ and $\log(1+\frac{\alpha}{\log x})$ is strictly 
decreasing and converges to zero, there exists $x_2$ such that 
$C>\log(1+\frac{\alpha}{\log x})$ for all $x\ge x_2$, hence by definition
of the number $e$ as an increasing limit
\begin{equation*}
\begin{split}
e^{C/\log x} &> 1+\frac{C}{\log x} > 1+\frac{\log(1+\frac{\alpha}{\log x})}{\log x}\\
&=\frac{\log x+\log(1+\frac{\alpha}{\log x})}{\log x}\\
&=\frac{\log(x(1+\frac{\alpha}{\log x}))}{\log x}
=\log_x\Big(x\big(1+\frac{\alpha}{\log x}\big)\Big)\,.
\end{split}
\end{equation*}
Therefore $x^{e^{C/\log x}}>x(1+\frac{\alpha}{\log x})$ and for 
a suitable $\alpha$ we can ensure that $x^{e^{C/\log x}}>\theta(x)$
for every $x\ge x_2$, we can choose for instance $\alpha = 1/2$ 
(see \cite[Theorem 4(3.15)]{Ros-Sch}). In consequence, taking
logarithms twice at both sides, we first get $e^{C/\log x}\log x>\log \theta(x)$, and then
$\log\Big(\frac{\log x}{\log\theta(x)}\Big)+\frac{C}{\log x}>0$ for every $x\ge x_2$.
This implies that there exists some $\varepsilon>0$ such that 
$$
\log\Big(\frac{\log x}{\log\theta(x)}\Big)+\frac{C}{\log x}-\varepsilon>0\,.
$$
Now we take into account computations of \eqref{gamma}, equation
\eqref{zeta}, \eqref{part} and add the Euler constant 
$\gamma$ to both sides of this inequality 
for $x\ge x_3=\max\{x_0,x_1,x_2\}$
\begin{equation*}
\begin{split}
\gamma &< \gamma + \log\Big(\frac{\log x}{\log\theta(x)}\Big)+\frac{C}{\log x}-\varepsilon\\
& < \gamma - b + b + \log\Big(\frac{\log x}{\log\theta(x)}\Big)+O(1/\log x)-\varepsilon\\
& = \sum_{n=2}^\infty \frac{P(n)}{n} -\varepsilon 
- \log\log\theta(x) + \sum_{p\le x}\frac1{p}\\
& < \sum_{n=2}^\infty \frac{P_x(n)}{n} 
- \log\log\theta(x) + \sum_{p\le x}\frac1{p}
= \log\Big(\frac{N_x}{\varphi_x\log\log N_x}\Big)\,.
\end{split}
\end{equation*}
Therefore Nicolas inequality holds for every $x\ge x_3$, which is not
exactly (T), but it is in contradiction with (F), and as by \cite[Theorem 2]{Nic} 
there are only two possibilities, (T) must hold.\\

{\it Step 3}: (II).
Now $C\le0$ and by hypothesis there is $x_0$ such that
$O(1/\log x)<\frac{C}{\log x}\le 0$ for every $x\ge x_0$, moreover we know
there is some $x_4\ge x_0$ such that $\theta(x_4)>x_4$
(see, for instance, the first page in \cite{Pla-Tru}), join
these facts with reasonings in previous steps, we obtain 
\begin{equation*}
\begin{split}
\gamma &\ge \gamma + \log\Big(\frac{\log x_4}{\log\theta(x_4)}\Big)+\frac{C}{\log x_4}\\
&\ge \sum_{n=2}^\infty \frac{P(n)}{n} - \log\log\theta(x_4) + \sum_{p\le x_4}\frac1{p}\\
&\ge \sum_{n=1}^\infty \frac{P_{x_4}(n)}{n} - \log\log\theta(x_4)\\
&= \log\Big(\frac{N_{x_4}}{\varphi_{x_4}\log\log N_{x_4}}\Big)\,,
\end{split}
\end{equation*}
hence Nicolas inequality does not hold for $x_4$ and we can apply the same argument to 
every integer $x\ge x_0$ such that $\theta(x)>x$ which is an infinite set of integers, i.e. assertion 
(F) holds.
\end{proof}

\begin{remark}\label{Nar}
Notice that conditions of theorem above are too restrictive, in the sense that condition of the 
first part never holds taking into account equation \cite[(2.3)]{Ros-Sch}, that is
$$
\sum_{p\le x} \frac1{p} = \log\log x + b + O(\exp(-\alpha\sqrt{\log x}))\,,
$$
for some positive constant $\alpha$. In fact 
\end{remark}

Observe that the failure of Nicolas inequality depends directly
from the change of sign of $\theta(x)-x$. From a revision of proof above we can deduce a
\begin{corollary}
Let $x\ge 2$, and let the Mertens formula
\begin{equation}\label{mertens}
\sum_{p\le x} \frac1{p} = \log\log x + b + O(1/\log x)\,.
\end{equation}
(I) Assume that $\theta(x)> x$. If $O(1/\log x)< 0$, then \eqref{Nicineq} does not hold for $x$.\\
(II) Assume that $\theta(x)\le x$. If $O(1/\log x)\ge 0$, then \eqref{Nicineq} holds for $x$.
\end{corollary}

An alternative to theorem above is

\begin{theorem}
Let $x> 2$, and let the Mertens formula
\begin{equation}\label{mertens}
\sum_{p\le x} \frac1{p} = \log\log x + b + O(1/\log x)\,.
\end{equation}
We define
$$
B(x):=\sum_{p> x} \log\big(1-\dfrac1{p}\big)+\dfrac1{p}\,,
$$
(I) Assume that $\theta(x)>x$. If \eqref{Nicineq} holds for $x$, then $O(1/\log x)\ge -B(x)$.\\
(II) Assume that $\theta(x)\le x$. If \eqref{Nicineq} does not hold for $x$, then 
$O(1/\log x) \le -B(x)\,.$
\end{theorem}
\begin{proof} We deduce (I) and (II) from a the following reasonings.
By rewriting proof of \cite[Theorem 8.8(d)]{Niv} we see that
$$
\prod_{p\le x}\big(1-\frac1{p}\big) = \frac{e^{-\gamma}}{\log x} \exp(-(O(1/\log x)+B(x)))\,,
$$
where $O(1/\log x)$ is in \eqref{mertens}, is the error term of Mertens Theorem and 
$$
B(x)=\sum_{p> x} \log\big(1-\dfrac1{p}\big)+\dfrac1{p}\,,
$$
it is an exercise to show that $B(x)<0$.
In our notation this can be written as
$$
C_x:=\frac{N_x}{\varphi_x\log x} = e^{\gamma} \exp(O(1/\log x)+B(x)) =  e^{\gamma} + \varepsilon_x\,,
$$
for some $\varepsilon_x\in (-1,1)$ for every $x> 2$.
We assumed that $x> 2$ is fixed.\\

(I) We first assume that $\theta_x=\log N_x > x$, 
then taking logarithms at both sides and multiplying by $C_k$ we obtain
$$
C_x\log\log N_x > \frac{N_x}{\varphi_x}\,,\quad\mbox{and hence }\quad
(e^\gamma + \varepsilon_x) \log\log N_x > \frac{N_x}{\varphi_x}\,.
$$
In consequence, by assuming that \eqref{Nicineq} holds for x
$$
\varepsilon_x\log\log N_x > \frac{N_x}{\varphi_x} - e^\gamma \log\log N_x > 0\,,
$$
since \eqref{Nicineq} holds. So $\varepsilon_x > 0$, thus $C_x > e^\gamma$. From defintion
of $C_x$ above we deduce that $\exp(O(1/\log x)+B(x))>1$, and so
$O(1/\log x)>-B(x)>0$.  In fact, if we assume that $O(1/\log x)\le 0$, we get 
$$
\log\Big(\frac{\log x}{\log\theta(x)}\Big)+O(1/\log x) \le 0\,,
$$
and following the arguments of Step 3 in the proof of previous theorem we obtain
$$
\gamma \ge \gamma + \log\Big(\frac{\log x}{\log\theta(x)}\Big) + O(1/\log x)
= \log\Big(\frac{N_{x}}{\varphi_{x}\log\log N_{x}}\Big)\,,
$$
hence Nicolas inequality does not hold for $x$.\\

(II) Now we assume that $\theta_x \le x$. By the same argumentation as above, if we assume
that \eqref{Nicineq} does not hold for $x$, we obtain in this case that $O(1/\log x)\le -B(x)$, where $-B(x)$
is positive.\\
\end{proof}

\begin{remark}
Observe that by definition of $B(x)$ and by Mertens Theorem the relation 
between $O(1/\log x)\lesseqgtr - B(x)$ is equivalent to
$$
\sum_{p\le x} \frac1{p} \lesseqgtr b + \log\log x -\sum_{p>x} \log\big(1-\dfrac1{p}\big)+\dfrac1{p}
=b + \log\log x +\sum_{p>x} \log\big(\dfrac{p}{p-1}\big)-\sum_{p>x}\dfrac1{p}\,.
$$
These inequalities could bring us to other discussions.\\
\end{remark}

In the following subsection we have used the quick guide of
Mertens Theorem given in \cite{Gol}.

\subsection{Discussion on the sign of the error term}
{\it Step 1}. Let start with an estimate for the sum of logarithms 
(see \cite[Lemma 8.2 and Theorem 8.8(a)]{Niv} and also \cite{Gol}).
$$
\sum_{n\le x} \log n = x \log x - x + O_1(\log x)\,,
$$
where do not mind the sign of $O_1(\log x)$ because in next step
we will divide it by $x$.\\

{\it Step 2}. We follow with a sum involving the Mangoldt function
(see \cite[Theorem 8.8(a)]{Niv}).
$$
\sum_{d\le x} \frac{\Lambda(d)}{d} = 
\log x + O_2(1)\,,
$$
where 
$$
O_2(1) = \frac{O_1(\log x)}{x} + \frac{\displaystyle\sum_{d\le x} \Big\{\dfrac{x}{d}\Big\}\Lambda(d)}{x}\,.
$$
On the one hand 
we argue by computation that the second term of this sum is
between $0.41$ and $0.43$ for sufficiently large $x$. This approximation
has been computed with the formula
\begin{equation}\label{aprox1}
0.41<\frac{\displaystyle\sum_{n=1}^{800}\Big(\displaystyle\sum_{j=1}^{[\pi(k)^{1/n}]}\Big(\frac{k}{p_j^n}
-\Big[\frac{k}{p_j^n}\Big]\Big)\log p_j\Big)}{k}<0.43\,,
\end{equation}

for $k\in(10^4,10^5)$ and a list of primes $\{p_j\}_j$. Alternatively we found estimates from above and 
from below for the second term. We can assume that $x$ is integer, hence the reminder of 
$\dfrac{x}{d}$ is $1$ for $d\in(x/2,x)$, so $\Big\{\dfrac{x}{d}\Big\}=\dfrac1{d}$ for 
$d\in(x/2,x)$, also recall that $\Big\{\dfrac{x}{x}\Big\}=0$. Analogously 
for $d\in(x/3,x/2)$, the reminder is $2$, so in this case $\Big\{\dfrac{x}{d}\Big\}=\dfrac{2}{d}$
and again $\Big\{\dfrac{x/2}{x/2}\Big\}=0$.
Taking into account \cite[Theorem 8.3]{Niv}, we can assert that 
$\psi(x/2)<3/4\psi(x)$ and $\psi(x/3)<1/2\psi(x)$ for a sufficient large integer $x$, thus
\begin{equation*}
\begin{split}
\sum_{d\le x}\Big\{\frac{x}{d}\Big\}\Lambda(d)
&=\sum_{d\le x/3}\Big\{\frac{x}{d}\Big\}\Lambda(d)
+\sum_{x/3< d \le x/2} \frac{2}{d} \Lambda(d)
+\sum_{x/2< d \le x} \frac1{d} \Lambda(d)\\
&<\sum_{d\le x/2}\Lambda(d)
+\frac{6}{x+1} \sum_{x/3< d < x/2} \Lambda(d)
+\frac{2}{x+1} \sum_{x/2< d < x} \Lambda(d)\\
&< \psi(x/3) +\frac{6}{x+1}(\psi(x/2)-\psi(x/3))+\frac{2}{x+1}(\psi(x)-\psi(x/2))\\
&<\frac1{2} \psi(x) + \frac{6}{x+1}\psi(x)+ \frac{2}{x+1}\psi(x)=\frac{x+17}{2x+2}\psi(x).
\end{split}
\end{equation*}

In consequence, taking into account that $ax<\psi(x)<cx$ for
$a=0.780355...$ and $c=3/2 a $, for a sufficiently large $x$ (see \cite[Theorem 8.3]{Niv}), we deduce that 
$$
0< \frac{\displaystyle\sum_{d\le x}
 \Big\{\dfrac{x}{d}\Big\}\Lambda(d)}{x}<0.58526625\,,
$$
where $0.58526625<\dfrac{3}{4}a$. On the other hand, as $O_1(\log x)/x$ converges to zero,
there is a $x_2$ such that $-1<O_2(1)<-0.415=0.585-1$ for every $x\ge x_2$.\\

{\it Step 3}. Next we deal with a sum over primes (see \cite[Theorem 8.8(b)]{Niv}).
$$
\sum_{p\le x} \frac{\log p}{p} = \log x + O_3(1)\,,
$$
where
$$
O_3(1) = O_2(1) - \sum_{n=2}^\infty\sum_{p\le x^{1/n}}\frac{\log p}{p^n}\,. 
$$
We approximate the second term by means of a list of primes and the formula
\begin{equation}\label{aprox2}
\sum_{n=2}^{40}\left(\sum_{j=1}^{[\pi(10^5)^{1/n}]}\frac{\log p_j}{p_j^n}\right)=0.753365132271...\,.
\end{equation}
Let us obtain two coarse bounds. We know that the second term is bounded and, 
since $\log p/p^n$ is positive, it increases when the $x$ grows and the partial sums in $n$ grows. Moreover
it is bounded from above by $\displaystyle\sum_{j=2}^\infty \frac{\log j}{j(j-1)}$
(see proof of \cite[Theorem 8.8(b)]{Niv}), which is bounded by the integral
of $\dfrac{\log (x-1)}{(x-1)(x-2)}$ between $2$ and infinity. 
So, for sufficiently large $x$
\begin{equation*}
\begin{split}
0.633 & <\sum_{n=2}^{40}\Big(\sum_{j=1}^{[\pi(100)^{1/n}]} \frac{\log p_j}{p_j^n}\Big) 
<\sum_{n=2}^\infty\sum_{p\le x^{1/n}}\frac{\log p}{p^n} <
\sum_{j=2}^\infty \frac{\log j}{j(j-1)}\\
&< \int_2^\infty \frac{\log(j-1)}{(j-1)(j-2)}=\frac{\pi^2}{6}<1.645\,.
\end{split}
\end{equation*}
where $(p_j)_j$ denotes the sequence of primes.
This implies that there exists $x_3$ such that, for every $x\ge x_3$:
\begin{equation}\label{O3}
-1.645 < O_3(1) < -0.048 < 0.585-0.633
\end{equation}

{\it Step 4}. Finally we study the error term of the sum of reciprocal of 
primes (see \cite[Theorem 8.8(d)]{Niv}).
We define $R(x)=\displaystyle\sum_{p\le x} \dfrac{\log p}{p} - \log x$, so thanks
to the estimate given in previous steps we get $R(x)=O_3(1)$. Let
$$
\sum_{p\le x} \frac1{p} = \log\log x + b + O(1/\log x)\,,
$$
where
\begin{equation}\label{bigo}
O(1/\log x)=\frac{R(x)}{\log x} - \int_x^\infty \frac{R(t)}{t(\log t)^2} \,dt\,.
\end{equation}
First observe that
$$
\frac{r}{\log x} - \int_x^\infty \frac{s}{t(\log t)^2}=\frac{r-s}{\log x}\,.
$$
Having in mind the above estimate for $O_3(1)=R(x)$, we assume that
$r,s\in(-1.645,-0.048)$ and obtain that
$$
\frac{-1.597}{\log x} < O(1/\log x) < \frac{1.597}{\log x}\,,
$$
for every $x\ge x_0=\max\{x_1,x_2,x_3\}$.\\

Estimates above are trivial, but we have seen that $R(x)$ is decreasing almost everywhere (a.e.) and 
negative for sufficiently large $x$. Let the expression of $b$ given 
in proof of \cite[Theorem 427]{Har-Wri} (see also \cite[Propostion]{Gol}), that is 
\begin{equation*}
\begin{split}
b & =\int_2^\infty \frac{R(t)}{t(\log t)^2}\,dt + 1 - \log\log 2\\
&=\int_2^x \frac{R(t)}{t(\log t)^2}\,dt + 1 - \log\log 2 + \int_x^\infty \frac{R(t)}{t(\log t)^2}\,dt\,.
\end{split}
\end{equation*}
Applying Stieltjes integration by parts a.e., we obtain
$$
\int_2^x \frac{R(t)}{t(\log t)^2}\,dt 
= \frac{R(2)}{\log 2} - \frac{R(x)}{\log x} + \int_2^x \frac{R'(t)}{\log t}\,dt 
= -1 - \frac{R(x)}{\log x} + \sum_{p\le x} \frac1{p} + \int_2^x \frac{R'(t)}{\log t}\,dt\,, 
$$
Moreover $\displaystyle\int_2^x \frac{R(t)}{t(\log t)^2}\,dt$ is negative and strictly 
decreasing that converges to $b+\log\log 2-1$, hence 
$$
\int_2^x \frac{R(t)}{t(\log t)^2}\,dt  > b+\log\log 2-1\,,
$$ 
for sufficiently large $x$. Notice that $R(x)$ is derivable and integrable (a.e.) 
and so is its derivative, moreover $R(2)/\log 2 = -1/2$.
Therefore
\begin{equation*}
\begin{split}
-\int_x^\infty \frac{R(t)}{t(\log t)^2}\,dt & = - b + 1 - \log\log 2 + \int_2^x \frac{R(t)}{t(\log t)^2}\,dt\\
& = - b + 1 - \log\log 2 - 1 - \frac{R(x)}{\log x} + \sum_{p\le x} \frac1{p} + \int_2^x \frac{R'(t)}{\log t}\,dt\,.
\end{split}
\end{equation*}
With these things in hand we see that
$$
O(1/\log x)=\frac{R(x)}{\log x} - \int_x^\infty \frac{R(t)}{t(\log t)^2}\,dt =
- b - \log\log 2 + \sum_{p\le x} \frac1{p} + \int_2^x \frac{R'(t)}{\log t}\,dt \,.
$$

\subsection*{Conclusion}
Observe that $O(1/\log x)$ converges to zero and, from the last equation, it is the difference between an step function
with discontinuities on primes and an strictly monotone function. We already saw this kind of
functions in \cite{Pla-Tru} (for instance), namely $\pi(x)-\operatorname{li}(x)$, $\psi(x)-x$ and $\theta(x)-x$, 
and we know they change of sign infinitely asymptoticaly. So it seems reasonable to think that $O(1/\log x)$
also change of sign intinitely for sufficiently large $x$, but this is another issue we do not 
deal here (see \cite[Chapter V]{Ing} for more information on these irregularities).

% ------------------------------------------------------------------------

\subsection*{Acknowledgment}
We would like to thank Professors 
Wladyslaw Narkiewicz, particularly for Remark \ref{Nar}, and Enrique A. S\'anchez P\'erez for their comments and suggestions.

% ------------------------------------------------------------------------
\end{document}